\documentclass[a4paper,11pt]{article}

\usepackage{mathpazo}

\usepackage{mathrsfs}
\usepackage{amssymb,amsmath,amsthm,url}
\usepackage{verbatim}
\usepackage[small]{titlesec}

\theoremstyle{plain} 
\newtheorem{thm}{Theorem}
\newtheorem*{corollary}{Corollary}
\newtheorem{lemma}{Lemma}[section]
\newtheorem{prop}[lemma]{Proposition}

\theoremstyle{definition} 
\newtheorem{problem}{Problem}
\newtheorem{definition}[lemma]{Definition}
\newtheorem*{acknowledgement}{Acknowledgements}
 
\newcommand{\LEFT}{(1)}
\newcommand{\RIGHT}{(2)}

\newcommand{\op}{\mbox{\scriptsize op}}

\newcommand{\IG}{\textsf{IG}}
\newcommand{\RIG}{\textsf{RIG}}
\newcommand{\MaxSbgp}{\textsf{H}}

\newcommand{\st}{\::\:}

\newcommand{\im}{\textup{im}}

\newcommand{\GR}{\mathscr{R}}
\newcommand{\GL}{\mathscr{L}}
\newcommand{\GH}{\mathscr{H}}
\newcommand{\GD}{\mathscr{D}}

\title{On Maximal Subgroups of Free Idempotent Generated Semigroups}

\author{R. Gray\thanks{Supported by an an EPSRC Postdoctoral Fellowship, and partially supported by FCT and FEDER, project POCTI-ISFL-1-143 of Centro de \'{A}lgebra da Universidade de Lisboa, and by the project PTDC/MAT/69514/2006.}, 
N. Ruskuc}

\begin{document}

\maketitle

\begin{abstract}
We prove the following results:
(1) Every group is a maximal subgroup of some free idempotent generated semigroup.
(2)
Every finitely presented group is a maximal subgroup of some free idempotent generated semigroup arising from a finite semigroup.
(3)
Every group is a maximal subgroup of some free regular idempotent generated semigroup.
(4)
Every finite group is a maximal subgroup of some free regular idempotent generated semigroup arising from a finite regular semigroup.
As a technical prerequisite for these results we establish a general presentation for the maximal subgroups based on a Reidemeister--Schreier type rewriting.
\smallskip

\end{abstract}

\section{Introduction and summary of results} 
\label{sec1}

Let $S$ be a semigroup, and let $E=E(S)$ be the set of idempotents of $S$.
The \emph{free idempotent generated semigroup} on $E$ is defined by the following presentation:
\begin{equation}
\label{eq1}
\IG(E)=\langle E\:|\: e\cdot f= ef\ (e,f\in E,\ \{ e,f\}\cap\{ef,fe\}\neq\emptyset)\rangle. 
\end{equation}
(It is an easy exercise to show that if, say, $fe\in\{e,f\}$ then $ef\in E$.
In the defining relation $e\cdot f=ef$ the left hand side is a word of length $2$, and $ef$ is the product of $e$ and $f$ in $S$, i.e. a word of length $1$.)
These semigroups arose in \cite{nambooripad79,easdown85},  where abstract characterisations of the sets of idempotents of semigroups via structures called \emph{biordered sets} was undertaken.
We will not need the formal definition of a biordered set here, but an interested reader may consult  \cite{higgins92}
for an accessible introduction.

The semigroup $\IG(E)$ has the following properties:

\renewcommand{\theenumi}{\textsf{(IG\arabic{enumi})}}
\renewcommand{\labelenumi}{\theenumi}

\begin{enumerate}
\item
\label{IG1}
There exists a natural homomorphism $\phi$ from $\IG(E)$ into the subsemigroup $S^\prime$ of $S$ generated by $E$.
\item
\label{IG2}
The restriction of $\phi$ to the set of idempotents of $\IG(E)$ is a bijection onto $E$ (and an isomorphism of biordered sets). Thus we may identify those two sets.
\item
\label{IG3}
$\phi$ maps the $\GR$-class (respectively $\GL$-class) of $e\in E$ onto the corresponding class of $e$ in $S^\prime$; this induces a bijection between the set of all $\GR$-classes (resp. $\GL$-classes) in the
$\GD$-class of $e$ in $\IG(E)$ and the corresponding set in $S^\prime$.
\item
\label{IG4}
The restriction of $\phi$ to the maximal subgroup of $\IG(E)$ containing $e\in E$ (i.e. to the $\GH$-class of $e$ in $\IG(E)$) is a homomorphism onto the maximal subgroup of $S^\prime$ containing $e$.
\end{enumerate}

The assertion \ref{IG1} is obvious; \ref{IG2} is proved in \cite{nambooripad79} and \cite{easdown85}; \ref{IG3} is a corollary of \cite{fitzgerald72}; \ref{IG4} follows from \ref{IG2}. 
The basics concerning Green's relations and their relationship with maximal subgroups will be reviewed in Section \ref{sec_prelims}. The maximal subgroup of a semigroup $S$ containing an idempotent $e$ will be denoted by $\MaxSbgp(S,e)$.

If $S$ is a regular semigroup, one also defines the \emph{free regular idempotent generated semigroup} $\RIG(E)$ on $E$ as follows.
The \emph{sandwich set} of a pair of idempotents $e,f\in E$ is defined as
$$
S(e,f)=\{ h\in E \st ehf=ef,\ fhe=h\}.
$$
The semigroup $\RIG(E)$ is then the homomorphic image of $\IG(E)$ obtained by adding the relations
$$
ehf=ef\ (e,f\in E,\ h\in S(e,f))
$$
to the presentation (\ref{eq1}). 
This semigroup also satisfies the properties \ref{IG1}--\ref{IG4}, and also:

\renewcommand{\theenumi}{\textsf{(RIG\arabic{enumi})}}
\renewcommand{\labelenumi}{\theenumi}

\begin{enumerate}
\setlength{\itemindent}{5mm}
\item
\label{RIG1}
$\RIG(E)$ is regular (\cite{nambooripad79}).
\item
\label{RIG2}
The natural homomorphism $\IG(E)\rightarrow \RIG(E)$ induces an isomorphism between the maximal subgroups of any $e\in E$ in $\IG(E)$ and $\RIG(E)$
(\cite[Theorem 3.6]{margolismeakin09}).
\end{enumerate}

Maximal subgroups of free idempotent generated semigroups have been of interest for some time.
Several papers \cite{mcelwee02,nambooripad80,pastijn77,pastijn80} established various sufficient conditions guaranteeing that all the maximal subgroups are free. Indeed, it was conjectured in \cite{mcelwee02} 
that this was always the case.
The first counterexample to this conjecture was given by Brittenham, Margolis and Meakin \cite{margolismeakin09}, where it was shown that the free abelian group $\mathbb{Z}\oplus\mathbb{Z}$ is the maximal subgroup of the free idempotent generated semigroup arising from a certain $72$-element semigroup. The same authors report further counterexamples to appear in \cite{margolismeakinip}, where they show that the multiplicative group $\mathbb{F}^\ast$ of a field $\mathbb{F}$ arises as a maximal subgroup of $\IG(E(M_3(\mathbb{F})))$, where $M_3(\mathbb{F})$ is the semigroup of all $3\times 3$ matrices over $\mathbb{F}$.
Related work concerning periodic elements in free idempotent generated semigroups may be found in \cite{easdownta}.

In this paper we prove:

\begin{thm}
\label{thm1}
Every group is a maximal subgroup of some free idempotent generated semigroup.
\end{thm}

\begin{thm}
\label{thm2}
Every finitely presented group is a maximal subgroup of some free idempotent generated semigroup arising from a finite semigroup.
\end{thm}

\begin{thm}
\label{thm3}
Every group is a maximal subgroup of some free regular idempotent generated semigroup.
\end{thm}

\begin{thm}
\label{thm4}
Every finite group is a maximal subgroup of some free regular idempotent generated semigroup arising from a finite regular semigroup.
\end{thm}

\begin{sloppypar}
We remark that Theorem \ref{thm2} provides a complete characterisation of groups appearing as maximal subgroups of free idempotent generated semigroups arising from finite semigroups. (And of course, trivially, Theorems \ref{thm1} and \ref{thm3} provide such characterisations with the finiteness assumption removed.)
Indeed, every maximal subgroup in a free idempotent generated semigroup arising from a finite semigroup must be finitely presented.
To see this, observe that in this case the presentation (\ref{eq1}) is finite, and also that the set
of $\mathcal{L}$-classes in the $\GD$-class of any $e\in E$ is finite by \ref{IG3}.
The assertion then follows from Proposition \ref{prop_finitelypresented} below.
By way of contrast, Theorem \ref{thm4} leaves us with the following unresolved question:
\emph{Is every finitely presented group a maximal subgroup of some free regular idempotent generated semigroup arising from a finite regular semigroup?}
\end{sloppypar}

Theorems \ref{thm1}--\ref{thm4} will be proved by means of two explicit constructions described in Sections \ref{sec3} and \ref{sec4}.
Preceding this, we prove a general presentation result based on the so called singular squares in Section \ref{sec_IGpres}, and introduce the features common to both constructions in Section \ref{sec2}.

\section{Preliminaries}
\label{sec_prelims}

In this section we review some basic definitions and facts concerning Green's relations, maximal subgroups, and a Reidemeister--Schreier type rewriting presentation for the latter. For a more systematic introduction to the basics of Semigroup Theory we refer the reader to any standard monograph such as \cite{higgins92,howie95}.

\subsection{Green's relations and maximal subgroups}
\label{subsec_grels}

Let $S$ be a semigroup. We use $S^1$ to denote the semigroup $S$ with an identity element $1 \not\in S$ adjoined to it. This notation will be extended to subsets of $S$, i.e. $X^1 = X \cup \{ 1 \}$. Green's relations, originally introduced in \cite{Green1951}, are equivalence relations which reflect the ideal structure of a semigroup. 
For $u,v \in S$ we define   
\[
u \GR v  \ \Leftrightarrow  \ uS^1 = vS^1, 
\quad 
u \GL v \  \Leftrightarrow \  S^1u = S^1v, \\
\]
\[
\GH = \GR \cap \GL, 
\quad
\GD = \GR \circ \GL = \GL \circ \GR.
\]
Each of these relations is an equivalence relation on $S$; their equivalence classes are called the
$\GR$-, $\GL$-, $\GH$- and $\GD$-classes, respectively. Moreover, $\GL$ is a right congruence and dually $\GR$ is a left congruence.  
The corresponding equivalence classes of an element $a \in S$ will be denoted by $R_a$, $L_a$, $H_a$ and $D_a$ respectively. 

Let $e$ be an idempotent of a semigroup $S$. The set $eSe$ is a submonoid of $S$ and is the largest submonoid whose identity is $e$. The group of units $G_e$ of $eSe$ (i.e. the group of elements of $eSe$ that have two-sided inverses with respect to $e$) is the largest subgroup of $S$ whose identity is $e$, and is called the maximal subgroup of $S$ containing $e$.  

We now list some fundamental facts about Green's relations and maximal subgroups, for proofs we refer the reader to 
\cite[Section~2]{howie95}.

\renewcommand{\theenumi}{\textsf{(G\arabic{enumi})}}
\renewcommand{\labelenumi}{\theenumi}

\begin{enumerate}
\item If $s,t \in S$ are such that $st \GR s$ then the mapping $\rho_t: x \mapsto xt$ is an $\GH$-class preserving bijection between the $\GL$-classes $L_s$ and $L_{st}$. 
\label{G1}
\item Furthermore, if $stu = s$, then the mappings $\rho_t$ and $\rho_u: L_{st} \rightarrow L_s$, $x \mapsto xu$, are mutually inverse bijections. 
\label{G2}
\item There are left/right dual statements to \ref{G1} and \ref{G2}. 
\label{G3}
\item The maximal subgroups of $S$ are precisely the $\GH$-classes that contain idempotents. 
\label{G4}
\item For any two $s, t \in S$ with $s \GD t$ we have $st \in R_s \cap L_t$ if and only if $R_t \cap L_s$ contains an idempotent. 
\label{G5}
\item If $s,t \in S$ are such that $st = s$ (resp. $ts=s$) then $xt = x$ (resp. $tx = x$) for all $x \in L_s$ (resp. $x \in R_s$). 
\label{G6}
\item In particular every idempotent is a left identity in its $\GR$-class and a right identity in its $\GL$-class. 
\label{G7}
\item If $s \in S$ is a regular element (meaning $sus = s$ for some $u \in S$) then all the elements in the $\GD$-class $D_s$ are regular (in which case we say that $D_s$ itself is regular). 
\label{G8}
\item In a regular $\GD$-class every $\GR$-class and every $\GL$-class contains an idempotent. 
\label{G9}
\item Let $e,f \in E(S)$. If $fe = e$ then $ef \in E(S)$ and $ef \GR e$. Dually, if $fe = f$ then $ef \in E(S)$ and $ef \GL f$. 
\label{G10}
\end{enumerate}

\subsection{A Reidemeister--Schreier type presentation for maximal subgroups}
\label{subsec_RSrewriting}

As a consequence of \ref{G1}, a semigroup acts on the $\GH$-classes in an $\GR$-class in a way similar to a group acting on the cosets of a subgroup. This analogy was exploited in \cite{ruskuc99} to obtain a presentation for an arbitrary maximal subgroup of a semigroup, closely resembling the classical Reidemeister--Schreier presentation in group theory. Let us now review this presentation.  

Let $S$ be a semigroup, and let $H$ be a maximal subgroup of $S$. Denote by $e$ the identity of $H$ (so $H = H_e$, the $\GH$-class of $e$). Let $R$ be the $\GR$-class of $e$, and let $H_j$ $(j \in J)$ be the $\GH$-classes contained in $R$. Without loss of generality assume $1 \in J$ and $H_1 = H$.    

The natural action of $S$ on itself by right multiplication induces an action on $R \cup \{ 0 \}$ (where $0$ is a new symbol intuitively meaning `undefined'). This action respects the $\GH$-classes in $\GR$ by \ref{G1}, and so can be represented by an action $(j,s) \mapsto j \cdot s$ of $S$ on $J \cup \{ 0 \}$ where
\[
j \cdot s = 
\begin{cases}
l & \mbox{if $j,l \in J$ and $H_j s = H_l$}, \\
0 & \mbox{otherwise}. 
\end{cases}
\] 
Suppose now that $S$ is generated by a set $A$. Let $r_j$ $(j \in J)$ be elements of $S^1$ such that 
\begin{equation}
\mbox{$H_1 r_j = H_j$ (or equivalently, $1 \cdot r_j = j$)}
\label{Rew1}
\end{equation}
for all $j \in J$.
By \ref{G1}, \ref{G2} there exist $r_j'$ $(j \in J)$ such that   
\begin{equation}
h r_j r_j' = h, \quad h' r_j' r_j = h' \quad (h \in H, \; h' \in H_j).
\label{Rew2}
\end{equation}
A (group) generating set for $H$ is given by
\begin{equation}
\{ 
e r_j a r_{j \cdot a}' :
j \in J, \; a \in A, \; j \cdot a \neq 0
\},
\label{Rew3}
\end{equation}
see \cite[Theorem~2.7]{ruskuc99}.

Now suppose further that $S$ is given by a presentation $\langle A | R \rangle$. Accordingly, suppose that the elements $e, r_j, r_j'$ $(j \in J)$ are given as words over $A$. Introduce a new alphabet
\begin{equation}
B = \{ 
[j,a] : j \in J, \; a \in A, \; j \cdot a \neq 0
\}
\label{Rew4}
\end{equation}
representing the generators \eqref{Rew3}. Denote by $\epsilon$ the empty word.
Define a \emph{rewriting mapping}
\[
\phi : \{  
(j,w) : j \in J, \; w \in A^*, \; j \cdot w \neq 0
\}
\rightarrow B^*
\]
inductively by
\begin{equation}
\phi(j,\epsilon) 	 = \epsilon, \quad\quad
 \phi(j,aw)  = [j,a] \phi(j \cdot a, w). 
\label{Rew5}
\end{equation}
Next suppose, without loss of generality, that the words $r_j$ $(j \in J)$ form a \emph{Schreier system}, i.e. every prefix of every $r_j$ is equal to some other $r_k$ (in particular this forces $r_1 = \epsilon$). A (group) presentation for $H$ is now given by
\begin{align}
H = \langle B \; | \; 
& 
[j,a] = 1
\ (j \in J, \; a \in A, \; j \cdot a \neq 0, \; r_{j \cdot a} = r_j a),
\label{Rew6} \\
& 
\phi(j,u) = \phi(j,v)
\ (j \in J, \; (u=v) \in R, \; j \cdot u \neq 0) \rangle,
\label{Rew7}
\end{align}
see \cite[Corollary~2.15]{ruskuc99}. Intuitively, this presentation for $H$ is obtained by rewriting the defining relations of $S$ with respect to the action on $\GH$-classes in $R$, and equating certain generators to $1$. An immediate consequence of this presentation is the following:
\begin{prop}\cite[Corollary~2.8]{ruskuc99}
\label{prop_finitelypresented}
Let $S$ be a semigroup and let $H$ be a maximal subgroup of $S$. If $S$ is finitely presented, and the number of $\GH$-classes in the $\GR$-class of $H$ is finite, then $H$ is finitely presented. 
\end{prop}

\section{Singular squares and a presentation for maximal subgroups in free idempotent generated semigroups}
\label{sec_IGpres}
\begin{sloppypar}
We are now going to employ the Reidemeister--Schreier type presentation described above to exhibit a `canonical' presentation for the maximal subgroup $H(\IG(E(S)),e)$ where $e \in E$. That this is a presentation for $H(\RIG(E(S)),e)$ ($S$ regular) was established by Nambooripad in \cite{nambooripad79} using different methods, and utilised in \cite{margolismeakin09,margolismeakinip}.
\end{sloppypar}

Let $S$ be an arbitrary semigroup, let $E = E(S)$, and form the free idempotent generated semigroup $\IG(E)$ as described in Section~\ref{sec1}. Clearly, without loss of generality we may assume that $S$ is generated by $E$ (i.e. $S = S'$ in the terminology of \ref{IG1}-\ref{IG4}).
Let $e_{1 1} \in E$ be arbitrary. We seek to obtain a presentation for the maximal subgroup $H = H(\IG(E),e_{1 1})$. A major significant point here is the following: the action of any generator $e \in E$ on the $\GH$-classes contained in the $\GR$-class of $e_{11}$ in $\IG(E)$ is equivalent to the action of $e$ on the $\GH$-classes contained in the $\GR$-class of $e_{11}$ in $S$. This follows from \ref{IG3}. Thus, in the presentation \eqref{Rew6}, \eqref{Rew7}, the relations to be rewritten come from the presentation \eqref{eq1} for $\IG(E)$, but the underlying action may be taken to come from $S$. 

Let us fix some notation. Let $D$ be the $\GD$-class of $e_{11}$, let $R_i$ $(i \in I)$ be the $\GR$-classes in $D$, let $L_j$ $(j \in J)$ be the $\GL$-classes in $\GD$, and let $H_{i j} = R_i \cap L_j$ for $i \in I, j \in J$. Next let $K = \{ (i,j) \in I \times J : \mbox{$H_{ij}$ is a group}  \}$, and let $e_{ij}$ be the identity of $H_{ij}$ for every $(i,j) \in K$. 
Focusing on the $\GR$-class $R = R_1$ and its $\GH$-classes $H_j = H_{1j}$ $(j \in J)$, let $r_j \in E^*$ be a Schreier system as described in Section~\ref{sec_prelims}. 
In fact, by the result of Fitzgerald \cite{fitzgerald72}, every element of $D$ can be expressed as a product of idempotents from $D$
of the form $e_{i_1j_1}e_{i_2j_2}\dots e_{i_nj_n}$ where $(i_{q+1},j_q)\in K$ for all $q=1,\dots,n-1$.
It is then easy to verify that we can choose our Schreier system $\{r_j\::\: j\in J\}$ to consist entirely of such products.

Keeping in mind the defining presentation \eqref{eq1} for $\IG(E)$ and the definition \eqref{Rew5} of the rewriting mapping $\phi$, the presentation \eqref{Rew6}, \eqref{Rew7} has generators 

\begin{equation}
B = 
\{   
[j,e] : j \in J, \; e \in E, \; j \cdot e \neq 0
\}                
\label{Rew8}
\end{equation}
and relations
\begin{align}
&[j,e] = 1 && (j \in J, \; e \in E, \; j \cdot e \neq 0, \; r_{j \cdot e} = r_j e),
\label{Rew9} \\
 \label{Rew10} 
& [t,e][t \cdot e,f] = [t,ef] &&
 (t \in J, \; e,f \in E, \; \{ e,f \} \cap \{ef, fe \} \neq \varnothing, \; t \cdot ef \neq 0).
\end{align}
Denote this presentation by $\mathcal{P}$.

In what follows we will first deduce some consequences of the relations \eqref{Rew9}, \eqref{Rew10},
leading to a new presentation, ultimately given in Theorem~\ref{lismaintheorem}. Then we will show that this new presentation actually implies all of the relations \eqref{Rew9}, \eqref{Rew10} and therefore it defines $H$. 

Let us partition the generators $B$ into two types:
\begin{align}
& B_1 = \{ [t,e] : e = e_{ij} \in D \ \& \ (i,t) \in K \},
\notag
\\ 
& B_2 = B \setminus B_1.
\notag
\end{align}
We remark that the condition $r_{j\cdot e}=r_j e$ appearing in (\ref{Rew9}) is an equality of two words over $E^\ast$,
and so,  because of the way in which we chose the Schreier system $\{r_j\::\: j\in J\}$,
all the generators appearing in (\ref{Rew9}) belong to $B_1$.

Let us, for every $i \in I$, fix $j(i) \in J$ such that $(i,j(i)) \in K$, which is possible by \ref{G9}.

\begin{lemma}
\label{lislemma1}
If $[t,e_{ij}]\in B_1$ then the relation
\begin{equation}
\label{liseq6}
[t,e_{ij}]=[j(i),e_{it}]^{-1}[j(i),e_{ij}]
\end{equation}
is a consequence of the presentation $\mathcal{P}$.
\end{lemma}

\begin{proof}
From $[t,e_{ij}]\in B_1$ it follows that $(i,t)\in K$.
The idempotents $e_{it}$ and $e_{ij}$ are $\GR$-related, and so the relation $e_{it}e_{ij}=e_{ij}$
is in presentation (\ref{eq1}) for $\IG(E)$. Also $(i,j(i))\in K$, and so
$e_{i,j(i)}e_{ij}=e_{ij}$, implying $j(i)\cdot e_{ij}=j\neq 0$ by \ref{G5}.
Hence the relation $[j(i),e_{it}][t,e_{ij}]=[j(i),e_{ij}]$ is in (\ref{Rew10}),
and it clearly implies the desired relation.
\end{proof}

\begin{lemma}
\label{lislemma1.5}
For any $(i,t)\in K$ the relation
\begin{equation}
\label{liseq7}
[t,e_{it}]=1
\end{equation}
is a consequence of \eqref{liseq6} (and hence of $\mathcal{P}$).
\end{lemma}

\begin{proof}
Put $j=t$ in Lemma \ref{lislemma1}.
\end{proof}

\renewcommand{\theenumi}{\textup{(\roman{enumi})}}
\renewcommand{\labelenumi}{\theenumi}

\begin{lemma}
\label{lislemma2}
Suppose $[t,e]\in B_2$ with $t\cdot e=j\in J$.
\begin{enumerate}
\item
\label{lislemma2item1}
There exists $i\in I$ such that $(i,t),(i,j)\in K$ and $ex=x$ for all $x\in R_i$.
\item
\label{lislemma2item2}
For any such $i$ the relation
\begin{equation}
\label{liseq8}
[t,e]=[t,e_{ij}]
\end{equation}
is a consequence of $\mathcal{P}$.
\end{enumerate}
\end{lemma}

\begin{proof}
\ref{lislemma2item1}
Pick any idempotent $e_{kj}\in L_j$ (it exists by \ref{G9}).
From $t\cdot e=j$ it follows that $L_te=L_j$, and hence $e$ stabilises  $L_j$ pointwise on the right by
\ref{G6}.
In particular, $e_{kj}e=e_{kj}$, and it follows by \ref{G10} that $ee_{kj}$ is also an idempotent in $L_j$,
say $e_{ij}$. Note that necessarily $(i,j)\in K$.
Now we have $ee_{ij}=eee_{kj}=ee_{kj}=e_{ij}$, and hence $ex=x$ for all $x\in R_i$ by \ref{G6}.
Also
$$
H_{1t}e_{ij}=H_{1t}ee_{ij}=H_{1j}e_{ij}=H_{1j},
$$
by \ref{G7}, and so \ref{G5} implies that $(i,t)\in K$.

\ref{lislemma2item2}
From $ee_{it}=e_{it}$ it follows that $e_{it}e$ is an idempotent in $R_i$.
Also, $e_{it}e\in L_t e=L_{t\cdot e}=L_j$, and hence $e_{it}e=e_{ij}$, a defining relation of
$\IG(E)$. Its rewritten form is $[t,e_{it}][t\cdot e_{it},e]=[t,e_{ij}]$ which belongs to (\ref{Rew10}).
By Lemma \ref{lislemma1.5} the relation $[t,e_{it}]=1$ is a consequence of $\mathcal{P}$.
Also, clearly $t\cdot e_{it}=t$, and the result follows.
\end{proof}

Relations in Lemma \ref{lislemma2} allow us to replace any generator from $B_2$ by a generator from $B_1$.
We observed earlier that no generator from $B_2$ appears in (\ref{Rew9}), so these relations remain unaltered by this substitution.
Relations in Lemma \ref{lislemma1} allow us to express all the generators from $B_1$ in terms of the subset
$$
\{ [j(i),e_{ij}] \::\: (i,j)\in K\}.
$$
Let us introduce a new notation for these generators:
$$
f_{ij}=[j(i),e_{ij}]\ ((i,j)\in K).
$$
The relations from Lemma \ref{lislemma1} become:
\begin{equation}
\label{liseq1}
[t,e_{ij}]=f_{it}^{-1} f_{ij}\ ((i,t),(i,j)\in K);
\end{equation}
they tell us how to replace any generator from $B_1$ in terms of the $f_{ij}$.
In order to do the same for generators from $B_2$ we first need to use Lemma \ref{lislemma2} to find a corresponding generator from $B_1$, and then use (\ref{liseq1}).

We now establish an important set of relations which the generators $f_{ij}$ satisfy. 

\begin{sloppypar}
\begin{definition}
\label{lisdef1}
A quadruple $(i,k;j,l)\in I\times I\times J\times J$ is a \emph{square}
if $(i,j),(i,l),(k,j),(k,l)\in K$.
It is a \emph{singular square} if, in addition, there exists an idempotent $e\in E$ such that
one of the following dual conditions holds:
\begin{eqnarray}
\label{liseq2}
&&ee_{ij}=e_{ij},\ ee_{kj}=e_{kj},\ e_{ij}e=e_{il},\ e_{kj}e=e_{kl},\ \mbox{or}\\
\label{liseq3}
&&e_{ij}e=e_{ij},\ e_{il}e=e_{il},\ ee_{ij}=e_{kj},\ ee_{il}=e_{kl}.
\end{eqnarray}
\end{definition}
\end{sloppypar}

We will say that $e$ \emph{singularises} the square.
Let $\Sigma_{LR}$ (respectively $\Sigma_{UD}$) be the set of all singular squares for 
which condition \eqref{liseq2}
(resp. \eqref{liseq3}) holds, and
let $\Sigma=\Sigma_{LR}\cup\Sigma_{UD}$, the set of all singular squares.
We call the members of $\Sigma_{LR}$ the \emph{left-right} singular squares, and those of $\Sigma_{UD}$ the \emph{up-down} singular squares.

\begin{lemma}
\label{lislemma3}
For every singular square $(i,k;j,l)\in\Sigma$ the relation
\begin{equation}
\label{liseq9}
f_{ij}^{-1}f_{il}=f_{kj}^{-1}f_{kl}
\end{equation}
is a consequence of $\mathcal{P}$.
\end{lemma}

\begin{proof}
Suppose first that (\ref{liseq2}) holds. By Lemma \ref{lislemma2} \ref{lislemma2item2}
and (\ref{liseq1}) we have
$$
f_{ij}^{-1}f_{il}=[j,e_{il}]=[j,e]=[j,e_{kl}]=f_{kj}^{-1}f_{kl}
$$
as a consequence of $\mathcal{P}$.

Consider now the case where (\ref{liseq3}) holds.
The relation $ee_{ij}=e_{kj}$ belongs to the presentation (\ref{eq1}), its rewritten form
is
\begin{equation}
\label{liseq5}
[l,e][l\cdot e,e_{ij}]=[l,e_{kj}],
\end{equation}
and it belongs to (\ref{Rew10}).
From $ee_{ij}=e_{kj}$ it follows that $ee_{kj}=e_{kj}$ and hence $ex=x$ for all $x\in R_k$.
From $e_{il}e=e_{il}$ we have $l\cdot e=l$.
Hence, by Lemma~\ref{lislemma2}~\ref{lislemma2item2} and \ref{lislemma1.5}, we have
$[l,e]=[l,e_{kl}]=1$. 
The relation (\ref{liseq5}) is now equivalent to $[l,e_{ij}]=[l,e_{kj}]$,
which, keeping in mind (\ref{liseq1}), is precisely the relation we sought to deduce.
\end{proof}

\begin{lemma}
\label{lislemma4}
Every relation
\begin{equation}
\label{liseq10}
[t,e][t\cdot e,f]=[t,ef]\ (e,f\in E,\ t\cdot ef\neq 0,\ \{e,f\}\cap\{ef,fe\}\neq \emptyset)
\end{equation}
from (\ref{Rew10}) is a consequence of relations (\ref{liseq6}) (equivalently, (\ref{liseq1})), (\ref{liseq8}), (\ref{liseq9}) from Lemmas \ref{lislemma1}, \ref{lislemma2}, \ref{lislemma3}.
\end{lemma}

\begin{proof}
We distinguish eight cases depending on whether each of the three generators belongs to $B_1$ or $B_2$.
\smallskip

\textit{Case 1: $[t,e],[t\cdot e,f],[t,ef]\in B_1$.}
From $[t,e]\in B_1$ it follows that $e=e_{ij}$, $(i,t)\in K$ and $t\cdot e=j$.
Similarly, from $[t\cdot e,f]=[j,f]\in B_1$ it follows that
$f=e_{kl}$ for some $k\in I$, and that $(k,j)\in K$.
From $(k,j)\in K$ and \ref{G5} we have $ef=e_{ij}e_{kl}\in H_{il}$, and so $(i,l)\in K$.
By \ref{G5} we have $fe\in H_{kj}$, and, since $fe$ is an idempotent, $fe=e_{kj}$.
Now $\{e,f\}\cap \{ef,fe\}$ implies that $i=k$ or $j=l$ (or both).
If $i=k$, using (\ref{liseq1}), we have
\begin{eqnarray*}
&&
[t,e]=[t,e_{ij}]=f_{it}^{-1} f_{ij},\\
&&
[t\cdot e,f]=[j,e_{il}]=f_{ij}^{-1} f_{il},\\
&&
[t,ef]=[t,e_{il}]=f_{it}^{-1}f_{il},
\end{eqnarray*}
and (\ref{liseq10}) follows readily.
Similarly, if $j=l$, we have $[t,e]=[t,ef]=f_{it}^{-1}f_{ij}$ and
$[t\cdot e,f]=[j,e_{kj}]=1$ (by (\ref{liseq7})), and (\ref{liseq10}) follows again.
\smallskip

\textit{Case 2: $[t,e],[t\cdot e,f]\in B_1$, $[t,ef]\in B_2$.}
In Case 1 we have only used the assumptions that $[t,e],[t\cdot e,f]\in B_1$,
and it then followed that $[t,ef]\in B_1$ as well.
Hence the present case cannot occur.
\smallskip

\textit{Case 3: $[t,e],[t,ef]\in B_1$, $[t\cdot e,f]\in B_2$.}
Utilise $[t,e],[t,ef]\in B_1$ to get:
$$
e=e_{ij},\ (i,t)\in K,\ t\cdot e=j,\ ef=e_{il},\ j\cdot f=l.
$$
From $eff=ef$ it follows that $fef$ is an idempotent $\GL$-related to $ef$, say $fef=e_{kl}$.
Now $fe_{kl}=e_{kl}$, and hence $fx=x$ for all $x\in R_k$.
Furthermore, notice that $H_{1j}e_{kl}=H_{1j}fef=H_{1l}ef=H_{1l}$; hence by \ref{G5} we have $(k,j)\in K$.
From Lemma \ref{lislemma2} \ref{lislemma2item2} we have $[j,f]=[j,e_{kl}]$.
Using (\ref{liseq1}) we have
\begin{eqnarray*}
&&
[t,e]=[t,e_{ij}]=f_{it}^{-1} f_{ij},\\
&&
[t\cdot e,f]=[j,f]=[j,e_{kl}]=f_{kj}^{-1} f_{kl},\\
&&
[t,ef]=[t,e_{il}]=f_{it}^{-1}f_{il}.
\end{eqnarray*}
Recall that we have $\{e,f\}\cap \{ef,fe\}\neq \emptyset$.
If $ef=e$ we have $j=l$ and (\ref{liseq10}) follows immediately.
The case $ef=f$ does not occur, as $(i,j)\in K$ would imply $[t\cdot e,f]=[j,f]\in B_1$.
If $fe=e$ we have $i=k$ and (\ref{liseq10}) follows.
Finally, the case $fe=f$ would imply $f=e_{kj}$ and $l=j$, and hence $[j,f]\in B_1$,
which is impossible.
\smallskip

\textit{Case 4: $[t,e]\in B_2$, $[t\cdot e,f],[t,ef]\in B_1$.}
Let $t\cdot e=j$.
From $[t\cdot e,f]\in B_1$ it follows that $f=e_{il}$ for some $l\in J$ and that $(i,j)\in K$.
Since $[t,ef]\in B_1$ we have $ef\in L_t ef=L_{t\cdot ef}=L_l$,
hence $ef=e_{kl}$ for some $k\in I$, and also $(k,t)\in K$.
From $ee_{il}=ef=e_{kl}$ it follows by \ref{G1}, \ref{G3} that $eH_{ij}=H_{kj}$.
Furthermore $e_{ij}e=e_{ij}$ and so $ee_{ij}$ is an idempotent by \ref{G10};
hence $(k,j)\in K$.
Note that $e$, $(k,t)$ and $(k,j)$ satisfy the condition from Lemma \ref{lislemma2} \ref{lislemma2item1}.
Hence Lemmas  \ref{lislemma2} \ref{lislemma2item2} and \ref{lislemma1} give
\begin{eqnarray*}
&&[t,e]=[t,e_{kj}]=f_{kt}^{-1} f_{kj},\\
&&[t\cdot e,f]=[j,e_{il}]=f_{ij}^{-1}f_{il},\\
&&[t,ef]=[t,e_{kl}]=f_{kt}^{-1}f_{kl}.
\end{eqnarray*}
As usual, we have $\{e,f\}\cap\{ef,fe\}\neq \emptyset$.
We cannot have $ef=e$ as it implies $[t,e]\in B_1$.
Neither can we have $fe=e$, as it implies $ef\GR e$ (by \ref{G10}), and again $[t,e]\in B_1$.
If $ef=f$ then $i=k$ and (\ref{liseq10}) follows immediately.
If $fe=f$ it follows that $l\cdot e=l$,  so that $(i,k;j,l)$ is a singular square (belonging to $\Sigma_{UD}$, singularised by $e$), and hence the relation $f_{ij}^{-1}f_{il}=f_{kj}^{-1}f_{kl}$ belongs to (\ref{liseq9}).
Now we have:
$$
[t,e][t\cdot e,f]=f_{kt}^{-1}f_{kj} f_{ij}^{-1}f_{il}
=
f_{kt}^{-1}f_{kj} f_{kj}^{-1}f_{kl}=f_{kt}^{-1} f_{kl}= [t,ef].
$$

\textit{Case 5: $[t,e],[t\cdot e,f]\in B_2$, $[t,ef]\in B_1$.}
Let $t\cdot e=j$ and $j\cdot f=l$.
From $ef\in B_1$ it follows that $ef\in L_t ef=L_l$ so that $ef=e_{il}$, and also $(i,t)\in K$.
Since $eef=ef$ it follows that $ex=x$ for all $x\in R_i$, and thus $(i,j)\in K$
and $e_{it}e=e_{ij}$ by \ref{G10}.
From $eff=ef$ and \ref{G10} it follows that $fef=fe_{il}=e_{kl}\in L_l$ for some $k\in I$.
Next note that $H_{1j}e_{kl}=H_{1j}fef=H_{1l}ef=H_{1l}$, since $ef\in L_l$ and \ref{G7}.
It follows by \ref{G5} that $(k,j)\in K$.
Regarding $ef$ and $fe$ the following cases cannot occur: $ef=e$ because it implies $[t,e]\in B_1$;
$ef=f$ because it implies $[t\cdot e,f]\in B_1$; $fe=f$ because it implies $f=fef=e_{kl}$ and so $[t\cdot e,f]\in B_1$.
In the remaining case $fe=e$ we have $e_{kl}=fef=ef=e_{il}$ so that $i=k$, and then the customary combination of (\ref{liseq8}) and (\ref{liseq1}) gives
$$
[t,e][t\cdot e,f] = [t,e_{ij}][j,e_{kl}]= f_{it}^{-1} f_{ij} f_{kj}^{-1} f_{kl}
= f_{it}^{-1} f_{ij} f_{ij}^{-1} f_{il}=f_{it}^{-1} f_{il}=[t,ef].
$$

\textit{Case 6: $[t,e],[t,ef]\in B_2$, $[t\cdot e,f]\in B_1$.}
Let $t\cdot e=j$.
Use $[t\cdot e,f]\in B_1$ to deduce $f=e_{il}$ for some $l\in J$, and also $(i,j)\in K$.
From $e_{ij}e=e_{ij}$ and \ref{G10} deduce $ee_{ij}=e_{kj}$ for some $k\in I$.
By \ref{G1}, \ref{G3}, we have $ef=ee_{il}\in H_{kl}$.
But $ef$ is an idempotent, and so $(k,l)\in K$ and $ef=e_{kl}$.
Next note $H_{1t}e_{kl}=H_{1t}ef=H_{1j}f=H_{1l}$, so that by \ref{G5} we have $(k,t)\in K$.
This in turn implies $[t,ef]\in B_1$, a contradiction.
\smallskip

\textit{Case 7: $[t,e]\in B_1$, $[t\cdot e,f],[t,ef]\in B_2$.}
From $[t,e]\in B_1$ deduce $e=e_{ij}$ and $(i,t)\in K$.
Let $j\cdot f=l\in J$.
Then $ef\in H_{ij}f=H_{il}$. But $ef\in E$, so $(i,l)\in K$ and $ef=e_{il}$.
But that implies $[t,ef]\in B_1$, a contradiction.
\smallskip

\textit{Case 8: $[t,e],[t\cdot e,f],[t,ef]\in B_2$.}
Let $t\cdot e=j$ and $j\cdot f=l$.
By Lemma \ref{lislemma2} there exists $i\in I$ such that
$efx=x$ for all $x\in R_i$, $(i,t),(i,l)\in K$ and $[t,ef]=[t,e_{il}]$.
Since $ef$ stabilises $R_i$ pointwise on the left, it follows that $e$ does so as well.
By \ref{G10} the element $e_{it}e\in H_{ij}$ is an idempotent, thus $(i,j)\in K$,
and $[t,e]=[t,e_{ij}]$ by Lemma \ref{lislemma2}.
Since $f$ stabilises $L_l$ pointwise on the right, it follows by \ref{G10} that $fe_{il}=e_{kl}$ for some $k\in I$. Next note that $H_{1j}e_{kl}=H_{1j}fe_{il}=H_{1l}e_{il}=H_{1l}$, and so
$(k,j)\in K$ by \ref{G5}. Also from $e_{kl}=fe_{il}$ we have that $f$ stabilises $R_k$ pointwise on the left, and $e_{kj}f=e_{kl}$, so that $[t\cdot e,f]=[j,f]=[j,e_{kl}]$ by Lemma \ref{lislemma2}. So now we have:
$$
[t,e]=f_{it}^{-1}f_{ij},\ [t\cdot e,f]=f_{kj}^{-1}f_{kl},\ [t,ef]=f_{it}^{-1}f_{il}.
$$
Let us also record that $fR_i=R_k$, and $eR_k=efR_i=R_i$.
Consider now the different possibilities for $ef$ and $fe$. 

\textit{8.1: $ef=e$.}
We have $j=j\cdot e=j\cdot ef=l$, so that $[t\cdot e,f]=1$, and (\ref{liseq10}) follows.

\textit{8.2: $ef=f$.}
From $R_i=efR_i=fR_i=R_k$ we have $i=k$, and (\ref{liseq10}) follows easily by cancellation of terms.

\textit{8.3: $fe=e$.}
We again have $R_i=eR_i=feR_i=fR_i=R_k$ and (\ref{liseq10}) follows as in 8.2.

\textit{8.4: $fe=f$.}
Now we have $l\cdot e=j\cdot fe=j\cdot f=l$. Hence $(k,i;j,l)\in\Sigma_{UD}$ (singularised by $e$),
and the relation $f_{kj}^{-1} f_{kl}=f_{ij}^{-1} f_{il}$ is in (\ref{liseq8}).
Now we have
$$
[t,e][t\cdot e,f]=f_{it}^{-1}f_{ij}f_{kj}^{-1}f_{kl}
=f_{it}^{-1}f_{ij}f_{ij}^{-1} f_{il}
=f_{it}^{-1} f_{il}= [t,ef].
$$

This exhausts all the cases and completes the proof of the lemma.
\end{proof}

Now we are going to use Lemmas \ref{lislemma1}--\ref{lislemma4} to transform the presentation (\ref{Rew8}), (\ref{Rew9}), (\ref{Rew10}) for $H$ to an equivalent presentation.
First use Lemmas \ref{lislemma1}, \ref{lislemma2}, \ref{lislemma3} to add all the relations
(\ref{liseq1}), (\ref{liseq8}), (\ref{liseq9}) to our presentation.
Then use Lemma \ref{lislemma4} to remove the redundant relations (\ref{Rew10}).

Next we would like to use relations (\ref{liseq8}) to eliminate all the generators from $B_2$.
However, note that there may be more than one relation (\ref{liseq8}) involving the same generator 
$[t,e]\in B_2$. Suppose we have two such relations $[t,e]=[t,e_{ij}]$ and $[t,e]=[t,e_{kj}]$,
and that we decide to use the former to eliminate $[t,e]$. This yields a new relation
$[t,e_{ij}]=[t,e_{kj}]$, i.e. $f_{it}^{-1}f_{ij}=f_{kt}^{-1}f_{kj}$.
Recall that here $i$ and $k$ both satisfy the condition of part \ref{lislemma2item1} of Lemma \ref{lislemma2}, which is easily seen to imply that $(i,k;t,j)\in\Sigma_{LR}$, and so the above relation belongs to (\ref{liseq9}).
Therefore we may indeed eliminate all the generators from $B_2$ and all the relations from (\ref{liseq8}).
Moreover, since no other remaining relations contain generators from $B_2$, they remain unchanged after this elimination.

Next we turn to eliminating the generators from $B_1$ other than the $f_{ij}$.
This is done using (\ref{liseq1}).
The only relations that remain from this group are those where $t=j(i)$, and they are clearly equivalent to $f_{i,j(i)}=1$ ($i\in I$).
Also, after this substitution,
every relation $[j,e_{il}]=1$ from \eqref{Rew9} becomes
$f_{ij}=f_{il}$. 

We have completed the proof of the following main theorem of this section, which essentially says that the maximal subgroups of free idempotent generated semigroups are defined by the relations arising from singular squares:

\begin{thm}
\label{lismaintheorem}
Let $S$ be a semigroup with a non-empty set of idempotents $E$,
let $\IG(E)$ be the corresponding free idempotent generated semigroup,
let $e_{11}\in E$ be arbitrary, and let $H$ be the maximal subgroup of $e_{11}$ in $\IG(E)$.
With the rest of notation as introduced throughout this section, a presentation for $H$ is given by
\begin{align}
\label{liseq25}
&\langle f_{ij}\ ((i,j)\in K) &&|&&
f_{ij}=f_{il} && ((i,j),(i,l)\in K,\ r_je_{il}=r_{j\cdot e_{il}}),
\\
\label{liseq26}
& && &&f_{i,j(i)}=1 && (i\in I),
\\
\label{liseq27}
& && &&f_{ij}^{-1}f_{il}=f_{kj}^{-1}f_{kl} && ((i,k;j,l)\in\Sigma) \rangle.
\end{align}
\end{thm}

\section{Outline of the method}
\label{sec2}

\subsection*{The environment semigroup \boldmath{$B_{I,J}$}}
There are two ways to compose mappings from a set $X$ into itself: from left to right and from right to left;
we denote the two resulting semigroups by $T_X$ and $T_X^{\op}$ respectively
(or $T_n$ and $T_n^{\op}$ if $X=\{1,\ldots,n\}$).
All the semigroups we are going to construct will be subsemigroups of some
$$
B_{I,J}=T_I^{\op}\times T_J.
$$
A typical element of $\beta\in B_{I,J}$ has the form $\beta=(\beta^{\LEFT},\beta^{\RIGHT})$.
To aid remembering the different orders in which compositions are formed we will write $\beta^{\LEFT}$
to the left of its argument, and $\beta^{\RIGHT}$ to the right.
The semigroup $B_{I,J}$ has a unique minimal ideal
$$
R_{I,J}=\{ \rho_{ij}=(\rho_i,\rho_j) \::\: i\in I,\ j\in J\},
$$
where
\begin{eqnarray*}
&& \rho_i :I\rightarrow I,\ x\mapsto i,\\
&&\rho_j : J\rightarrow J,\ x\mapsto j,
\end{eqnarray*}
are the constant maps. The multiplication in $R_{I,J}$ works as follows:
$$
\rho_{ij}\rho_{kl}=\rho_{il},
$$
i.e. $R_{I,J}$ is an $I\times J$ rectangular band.
The semigroup $B_{I,J}$ is in fact the \emph{translational hull} of $R_{I,J}$
(see \cite{petrich68}), an important background fact for our discussion,
even though it will not be explicitly used in any of the arguments.

\subsection*{Maximal subgroups in the minimal ideal}

In each of the two constructions to be described in Sections \ref{sec3}, \ref{sec4} we work with
a semigroup $S$ satisfying:

\renewcommand{\theenumi}{\textsf{(S\arabic{enumi})}}
\renewcommand{\labelenumi}{\theenumi}

\begin{enumerate}
\item
\label{S1}
$R_{I,J}\leq S\leq B_{I,J}$.
\end{enumerate}
In each instance we will want to determine the maximal subgroup $H$
of $\IG(E(S))$ containing $\rho_{11}$.

The fact that $\rho_{11}$ is in a rectangular band  minimal ideal of $S$ implies some further simplifications to the general presentation given in Theorem \ref{lismaintheorem}, which we outline here.
First note that now $K=I\times J$, and so our generating set is simply 
\begin{equation}
\label{eq11}
f_{ij} \ (i\in I,\ j\in J).
\end{equation}
Next, note that for every $i\in I$ the element $j(i)$ can be simply taken to be $1$, whereupon
relations (\ref{liseq26}) become $f_{i1}=1$ ($i\in I$).
The words $r_1=\epsilon$, $r_j=\rho_{1j}$ ($1\neq j\in J$) clearly form a Schreier system of representatives. With this system, relations (\ref{liseq25}) become $1=f_{11}=f_{1l}$ ($1\neq l\in J$). 
So the group $H=\MaxSbgp(\IG(E(S)),\rho_{11})$ is defined by the presentation
\begin{align}
\label{liseq28}
&\langle f_{ij}\ (i\in I,\ j\in J) &&|&&
f_{1j}=f_{i1}=1 && (i\in I,\ j\in J),\\
\label{liseq29}
& && && f_{ij}^{-1}f_{il}=f_{kj}^{-1} f_{kl} && ( (i,k;j,l)\in\Sigma))\rangle .
\end{align}

Since our semigroups consist of (pairs of) mappings and the $\GD$-class under consideration consists of (pairs of) constant mappings, the definition of singular squares can be conveniently recast as follows:
\begin{equation}
\label{eq13a}
\begin{split}
\Sigma &= \Sigma_{LR}\cup \Sigma_{UD}\\
\Sigma_{LR} &= \{ (i,k;j,l)\in I\times I\times J\times J\::\: \\
&\quad\quad (\exists \beta\in E(S))(\beta^{\LEFT}(i)=i\ \&\ \beta^{\LEFT}(k)=k\ \&\ j\beta^{\RIGHT}=l\beta^{\RIGHT}=j)\}
\\
\Sigma_{UD} &= \{ (i,k;j,l)\in I\times I\times J\times J\::\: \\
&\quad\quad (\exists \beta\in E(S))(\beta^{\LEFT}(i)=\beta^{\LEFT}(k)=i\ \&\  j\beta^{\RIGHT}=j\ \&\ l\beta^{\RIGHT}=l)\}.
\end{split}
\end{equation}

\begin{sloppypar}
To aid understanding of the forthcoming considerations, let us highlight the relations produced by certain distinguished types of singular squares.
For example, a square $(1,i;1,j)\in \Sigma$ yields the relation $f_{11}^{-1}f_{1j}=f_{i1}^{-1}f_{ij}$,
which, keeping in mind (\ref{liseq28}), is equivalent to $f_{ij}=1$.
A singular square of the form $(1,i;j,l)$ yields the relation $f_{ij}=f_{il}$,
while the square $(i,k;1,j)$ yields $f_{ij}=f_{kj}$.
Let us call these three types the \emph{corner}, \emph{flush top} and \emph{flush left} squares, respectively.
\end{sloppypar}

One further type of square will be utilised: Suppose that $(i,k;j,l)\in\Sigma$ and that we already know that
$f_{ij}=1$ (e.g. by virtue of a corner square involving $i$ and $j$). Then the relation (\ref{liseq29})
becomes $f_{kj}f_{il}=f_{kl}$. Let us call this a \emph{$3/4$ square}. All four different types of singular squares are illustrated in Figure \ref{fig1}.

\setlength{\unitlength}{1mm}

\begin{figure}
\begin{center}

\begin{tabular}{ccc}
\begin{picture}(40,35)(0,5)
\put(5,5){\framebox(35,30){}}
\put(7,15){\makebox(3,3){$1$}}
\put(7,30){\makebox(3,3){$1$}}
\put(27,15){\makebox(3,3){$f_{ij}$}}
\put(27,30){\makebox(3,3){$1$}}

\thicklines
\put(8.5,16.5){\circle{6}}
\put(8.5,31.5){\circle{6}}
\put(28.5,16.5){\circle{6}}
\put(28.5,31.5){\circle{6}}
\put(11.5,16.5){\line(1,0){14}}

\put(11.5,30.5){\line(1,0){14}}
\put(8.5,19.5){\line(0,1){9}}
\put(28.5,19.5){\line(0,1){9}}

\put(1,15){\makebox(3,3){$i$}}
\put(1,30){\makebox(3,3){$1$}}
\put(7,37){\makebox(3,3){$1$}}
\put(27,37){\makebox(3,3){$j$}}

\put(15,20){\makebox(10,5){$f_{ij}=1$}}
\end{picture}
&
\makebox[10mm]{}
&
\begin{picture}(50,35)(0,5)
\put(5,5){\framebox(45,30){}}
\put(7,15){\makebox(3,3){$1$}}
\put(7,30){\makebox(3,3){$1$}}
\put(22,15){\makebox(3,3){$f_{ij}$}}
\put(22,30){\makebox(3,3){$1$}}
\put(42,15){\makebox(3,3){$f_{il}$}}
\put(42,30){\makebox(3,3){$1$}}

\thicklines

\put(23.5,16.5){\circle{6}}
\put(23.5,31.5){\circle{6}}
\put(43.5,16.5){\circle{6}}
\put(43.5,31.5){\circle{6}}

\put(23.5,19.5){\line(0,1){9}}
\put(43.5,19.5){\line(0,1){9}}
\put(26.5,16.5){\line(1,0){14}}
\put(26.5,31.5){\line(1,0){14}}

\put(1,15){\makebox(3,3){$i$}}
\put(1,30){\makebox(3,3){$1$}}
\put(7,37){\makebox(3,3){$1$}}
\put(22,37){\makebox(3,3){$j$}}
\put(42,37){\makebox(3,3){$l$}}

\put(28,20){\makebox(10,5){$f_{ij}=f_{il}$}}
\end{picture}
\\
Corner&&Flush top
\\
&&
\\
\begin{picture}(40,45)(0,5)
\put(5,5){\framebox(35,40){}}
\put(7,30){\makebox(3,3){$1$}}
\put(7,40){\makebox(3,3){$1$}}
\put(27,30){\makebox(3,3){$f_{ij}$}}
\put(27,40){\makebox(3,3){$1$}}
\put(7,10){\makebox(3,3){$1$}}
\put(27,10){\makebox(3,3){$f_{kj}$}}

\thicklines

\put(8.5,31.5){\circle{6}}
\put(28.5,31.5){\circle{6}}
\put(8.5,11.5){\circle{6}}
\put(28.5,11.5){\circle{6}}

\put(11.5,31.5){\line(1,0){14}}
\put(11.5,11.5){\line(1,0){14}}
\put(8.5,14.5){\line(0,1){14}}
\put(28.5,14.5){\line(0,1){14}}

\put(1,30){\makebox(3,3){$i$}}
\put(1,10){\makebox(3,3){$k$}}
\put(1,40){\makebox(3,3){$1$}}
\put(7,47){\makebox(3,3){$1$}}
\put(27,47){\makebox(3,3){$j$}}

\put(14,20){\makebox(10,5){$f_{ij}=f_{kj}$}}
\end{picture}
&&
\begin{picture}(50,45)(0,5)
\put(5,5){\framebox(45,40){}}
\put(7,30){\makebox(3,3){$1$}}
\put(7,40){\makebox(3,3){$1$}}
\put(22,30){\makebox(3,3){$1$}}
\put(22,40){\makebox(3,3){$1$}}
\put(7,10){\makebox(3,3){$1$}}
\put(22,10){\makebox(3,3){$f_{kj}$}}
\put(42,10){\makebox(3,3){$f_{kl}$}}
\put(42,30){\makebox(3,3){$f_{il}$}}
\put(42,40){\makebox(3,3){$1$}}

\thicklines

\put(23.5,31.5){\circle{6}}
\put(23.5,11.5){\circle{6}}
\put(43.5,11.5){\circle{6}}
\put(43.5,31.5){\circle{6}}

\put(23.5,14.5){\line(0,1){14}}
\put(43.5,14.5){\line(0,1){14}}
\put(26.5,11.5){\line(1,0){14}}
\put(26.5,31.5){\line(1,0){14}}

\put(1,30){\makebox(3,3){$i$}}
\put(1,10){\makebox(3,3){$k$}}
\put(1,40){\makebox(3,3){$1$}}
\put(7,47){\makebox(3,3){$1$}}
\put(22,47){\makebox(3,3){$j$}}
\put(42,47){\makebox(3,3){$l$}}

\put(25,20){\makebox(17,5){$f_{kj}f_{il}=f_{kl}$}}
\end{picture}
\\
Flush left&&$3/4$
\end{tabular}
\caption{The four distinguished types of singular squares, and the relations they yield.}
\label{fig1}
\end{center}
\end{figure}
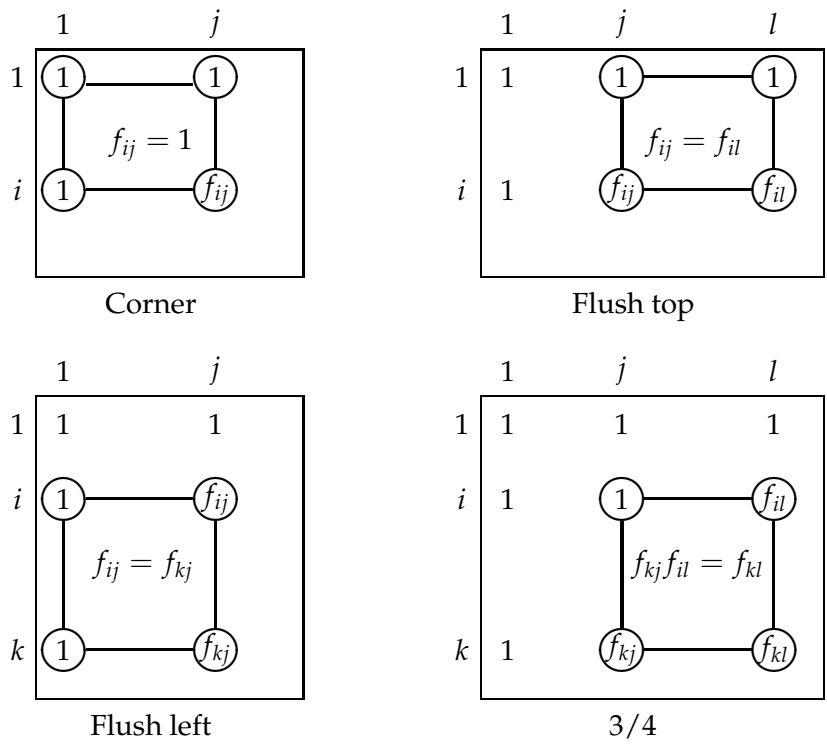

\subsection*{The method}
Let us outline the features and reasoning common to both forthcoming constructions.
We will start with a group $G$, which ultimately we want to realise as a maximal subgroup of the free idempotent generated semigroup arising from some semigroup $S$.
In each instance we will introduce an \emph{auxiliary matrix}
\begin{equation}
\label{eq13aa}
Y=(y_{ij})_{I\times J},
\end{equation}
with entries from a generating set for $G$.
Aided by this auxiliary matrix, we define a collection of idempotents from $B_{I,J}$,
which, together with $R_{I,J}$, generate $S$.
We then examine the presentation given in \eqref{liseq28}, \eqref{liseq29} and prove that it indeed defines our initial group $G$.
This we do by performing the following steps. 
First we consider certain idempotents from $S\setminus R_{I,J}$ such that the (corner and flush) squares singularised by them allow us to prove:

\renewcommand{\theenumi}{\textsf{(Rel\arabic{enumi})}}
\renewcommand{\labelenumi}{\theenumi}

\begin{enumerate}
\setlength{\itemindent}{3mm}
\item
\label{Rel1}
$f_{ij}=f_{kl}$ whenever $y_{ij}=y_{kl}$.
\end{enumerate}
At this stage we may identify each letter $f_{ij}$ with the corresponding entry $y_{ij}$ (considered also as a letter, and not as a group element).
Then we consider certain further idempotents such that:

\begin{enumerate}
\setlength{\itemindent}{3mm}
\setcounter{enumi}{1}
\item
\label{Rel2}
The ($3/4$) squares singularised by these idempotents yield (after the identification above) a set of relations which define $G$.
\end{enumerate}
The argument is completed by a `bookkeeping step', which ensures that we do get the group $G$ and not a proper homomorphic image:

\begin{enumerate}
\setlength{\itemindent}{3mm}
\setcounter{enumi}{2}
\item
\label{Rel3}
Check (by considering every idempotent and every square singularised by it) that every other relation from \eqref{liseq28}, \eqref{liseq29} holds in $G$.
\end{enumerate}

It in fact turns out that the auxiliary matrix $Y$ yields a natural Rees matrix representation
(see \cite[Section 3.2]{howie95})
$\mathcal{M}[H_{11};I,J;P]$ for the minimal ideal $D$ of $\IG(E)$, with the structure matrix $P=(p_{ji})_{J\times I}$ given by
$p_{ji}=y_{ij}^{-1}$.

\section{First construction}
\label{sec3}

In this section we present our first construction, which proves Theorems \ref{thm1} and \ref{thm2}.
In fact, we will present the construction in the finitely presented case (Theorem \ref{thm2}), and then
remark that obvious trivial changes adapt it to the general case (Theorem \ref{thm1}).

So let $G$ be any finitely presented group. It is well known that
$G$ can be defined by a presentation of the form
\begin{equation}
\label{eq14}
G=\langle a_1,\ldots,a_p\:|\: b_1c_1=d_1,\ldots, b_qc_q=d_q\rangle,
\end{equation}
where $b_r,c_r,d_r\in\{a_1,\ldots,a_p\}$ for all $r=1,\ldots,q$.
This essentially follows from the fact that a `long' relator $x_1x_2\ldots x_s$
($s>3$) can be replaced by two shorter relators $x_1x_2y^{-1}$ and $yx_3\ldots x_s$, at the expense
of introducing a new, redundant generator $y$.

Let us set
\begin{equation}
\label{eq15}
m=1+2q,\ n=1+p+2q,\ I=\{1,\ldots,m\},\ J=\{1,\ldots,n\}.
\end{equation}
The auxiliary matrix (see Section \ref{sec2}) is:
\begin{equation}
\label{eq15aa}
Y=(y_{ij})_{m\times n}=
\left( \begin{array}{cccccccccccc}
          1&1&1&\ldots&1&1&1&1&1&\ldots&1&1\\
          1&a_1&a_2&\ldots&a_p&1&c_1&1&1&\ldots&1&1\\
          1&a_1&a_2&\ldots&a_p&b_1&d_1&1&1&\ldots&1&1\\
          1&a_1&a_2&\ldots&a_p&1&1&1&c_2&\ldots&1&1\\
          1&a_1&a_2&\ldots&a_p&1&1&b_2&d_2&\ldots&1&1\\
          &\vdots&&&&&\vdots&&&\ddots&\vdots&\\
          1&a_1&a_2&\ldots&a_p&1&1&1&1&\ldots&1&c_q\\
          1&a_1&a_2&\ldots&a_p&1&1&1&1&\ldots&b_q&d_q
       \end{array}\right).
\end{equation}

Next we define a family of idempotents which, together with $R_{m,n}$, will 
be used to generate $S$:
\begin{equation}
\label{eq15a}
\begin{split}
&\sigma_u=(\sigma_u^{\LEFT},\sigma_u^{\RIGHT})\ (u=2,\ldots,m),\\
&\tau_u=(\tau_u^{\LEFT},\tau_u^{\RIGHT})\ (u=1,\ldots,q),
\end{split}
\end{equation}
where
\begin{eqnarray}
\label{eq15b}
\sigma_u^{\LEFT}(x)&=&\left\{ \begin{array}{ll} 1&\mbox{if } x=1\\
                                              u&\mbox{if } x\neq 1
                          \end{array}\right.
\\
\label{eq15c}
x\sigma_u^{\RIGHT} &=& \left\{ \begin{array}{ll} x&\mbox{if } 1\leq x\leq p+1\\
                                              r+1 &\mbox{if } x>p+1 \mbox{ and } y_{u,x}=a_r
                          \end{array}\right.
\end{eqnarray}
(with the convention that $a_0=1$), and
\begin{eqnarray}
\label{eq15d}
\tau_u^{\LEFT}(x)&=&\left\{ \begin{array}{ll} 2u+1&\mbox{if } x=2u+1\\
                                            2u&\mbox{if } x\neq 2u+1
                          \end{array}\right.
\\
\label{eq15e}
x\tau_u^{\RIGHT} &=& \left\{ \begin{array}{ll} p+2u&\mbox{if } x=p+2u,p+2u+1\\
                                            1 &\mbox{otherwise}.
                          \end{array}\right. 
\end{eqnarray}

Let us interrupt our exposition at this point in order to illustrate the construction thus far by means of a concrete example.
Let us take $G$ to be the Klein bottle group
$$
G=\langle a,b\:|\: a^{-1}ba=b^{-1}\rangle.
$$
By introducing a redundant generator $c=ba$ we obtain the following presentation
$$
G=\langle a,b,c \:|\: ba=c,\ cb=a\rangle,
$$
which has the form (\ref{eq14}).
Thus here $p=3$, $q=2$, $m=5$, $n=8$, and the auxiliary matrix is
$$
Y=\left( \begin{array}{cccccccc}
1&1&1&1&1&1&1&1\\
1&a&b&c&1&a&1&1\\
1&a&b&c&b&c&1&1\\
1&a&b&c&1&1&1&b\\
1&a&b&c&1&1&c&a
\end{array}\right).
$$

The only definition from (\ref{eq15b})--(\ref{eq15e}) that is perhaps not entirely straightforward is
(\ref{eq15c}).
So let us compute $\sigma_3^{\RIGHT}$. Its image is $\{1,2,3,4\}$, and it acts identically on it.
Each of the remaining columns is mapped into that image column which has the same 3rd entry.
The 3rd entries for columns $5$, $6$, $7$, $8$ are $b$, $c$, $1$, $1$ respectively, and so $5$, $6$, $7$, $8$ are mapped to $3$, $4$, $1$, $1$ respectively.

The complete list of our idempotents is:
\begin{alignat*}{5}
&\sigma_2^{\LEFT}&&= \left(\begin{array}{ccccc} 1&2&3&4&5\\1&2&2&2&2\end{array}\right) && \quad\quad
&&\sigma_2^{\RIGHT}&&= \left(\begin{array}{cccccccc} 1&2&3&4&5&6&7&8\\1&2&3&4&1&2&1&1\end{array}\right)
\\
&\sigma_3^{\LEFT}&&= \left(\begin{array}{ccccc} 1&2&3&4&5\\1&3&3&3&3\end{array}\right) && \quad 
&&\sigma_3^{\RIGHT}&&= \left(\begin{array}{cccccccc} 1&2&3&4&5&6&7&8\\1&2&3&4&3&4&1&1\end{array}\right)
\\
&\sigma_4^{\LEFT}&&= \left(\begin{array}{ccccc} 1&2&3&4&5\\1&4&4&4&4\end{array}\right) &&\quad 
&&\sigma_4^{\RIGHT}&&= \left(\begin{array}{cccccccc} 1&2&3&4&5&6&7&8\\1&2&3&4&1&1&1&3\end{array}\right)
\\
&\sigma_5^{\LEFT}&&= \left(\begin{array}{ccccc} 1&2&3&4&5\\1&5&5&5&5\end{array}\right) && \quad 
&&\sigma_5^{\RIGHT}&&= \left(\begin{array}{cccccccc} 1&2&3&4&5&6&7&8\\1&2&3&4&1&1&4&2\end{array}\right)
\\
&\tau_1^{\LEFT}&&= \left(\begin{array}{ccccc} 1&2&3&4&5\\2&2&3&2&2\end{array}\right) && \quad 
&&\tau_1^{\RIGHT}&&= \left(\begin{array}{cccccccc} 1&2&3&4&5&6&7&8\\1&1&1&1&5&5&1&1\end{array}\right)
\\
&\tau_2^{\LEFT}&&= \left(\begin{array}{ccccc} 1&2&3&4&5\\4&4&4&4&5\end{array}\right)  
&&\quad 
&&\tau_2^{\RIGHT}&&= \left(\begin{array}{cccccccc} 1&2&3&4&5&6&7&8\\1&1&1&1&1&1&7&7\end{array}\right).
\end{alignat*}

Returning to our general argument, let $S$ be the subsemigroup of $B_{m,n}$ generated
by the idempotents $\sigma_u$ ($u=2,\ldots,m$), $\tau_u$ ($u=1,\ldots,q$) and $R_{m,n}$.
We start working towards accomplishing \ref{Rel1} (see Section \ref{sec2}) by identifying the squares singularised by the idempotent $\sigma_u$ ($u=2,\ldots,m$) and the relations these squares yield.
Beginning with the left-right squares, notice that $\im(\sigma_u^{\LEFT})=\{1,u\}$, and so a typical square
singularised has the form $(1,u; r+1,l)$ where $0\leq r\leq p$, $l>p+1$ and $y_{u,l}=a_r$.
This is a flush top square, and yields the relation
\begin{equation}
\label{eq21}
f_{ul}=f_{u,r+1}\ (2\leq u\leq m,\ p+1<l<n,\ 0\leq r\leq p,\ y_{ul}=a_r).
\end{equation}
The up-down squares singularised by $\sigma_2$ are of the form
$(2,i; j,l)$ ($2< i\leq m,\ 1\leq j<l\leq p+1$).
For $j=1$ we obtain a flush left square $(2,i;1,l)$, which yields the relation
\begin{equation}
\label{eq23}
f_{il}=f_{2l}\ (2<i\leq m,\ 1<l\leq p+1).
\end{equation}
A general up-down square singularised by any $\sigma_u$ has the form
$(u,i; j,l)$ ($2\leq i\leq m,\ i\neq u,\ 1\leq j<l\leq p+1$) and the relation it yields is an easy consequence of
(\ref{eq23}):
$$
f_{uj}^{-1}f_{ul}=f_{2j}^{-1}f_{2l}=f_{ij}^{-1}f_{il}.
$$

We see that the squares singularised by $\sigma_u$ ($2\leq u\leq m$) allow us precisely to deduce \ref{Rel1}:
A typical generator $f_{ij}$ is equal to an appropriate $f_{i,j_1}$ (with $1\leq j_1 \leq p+1$) by (\ref{eq21}), which in turn is equal to its `canonical representative' $f_{2,j_1}$ by 
(\ref{eq23}).
So from now on we identify each generator $f_{ij}$ with $y_{ij}$.

Let us now examine the squares singularised by $\tau_u$ ($u=1,\ldots,q$).
The most important such square is $(2u,2u+1; p+2u,p+2u+1)$.
It is a left-right 3/4 square, since under our identification
$f_{2u,p+2u}=y_{2u,p+2u}=1$, and yields the relation
\begin{equation}
\label{eq25}
b_uc_u=d_u\ (u=1,\ldots,q).
\end{equation}
Every other left-right square for $\tau_u$ is flush left $(2u,2u+1; 1,j)$ ($j\neq 1,2u,2u+1$),
yielding either $a_{j-1}=a_{j-1}$ (for $j\leq p+1$) or $1=1$ (for $j>p+1$).
The up-down squares have the form $(2u,i;1,p+2u)$ ($i\neq 2u+1$), and yield the trivial relation $1=1$.

At this stage we have accomplished step \ref{Rel2}, and have precisely a presentation defining $G$.
Furthermore we have examined all the squares singularised by $\sigma_2,\ldots,\sigma_m,\tau_1,\ldots,\tau_q$.
So, to verify \ref{Rel3}, i.e. prove that no further relations are introduced, it is sufficient to prove that $S$ contains no further idempotents:  
\begin{equation}
\label{eq26}
E(S)=R_{m,n}\cup\{\sigma_2,\ldots,\sigma_m,\tau_1,\ldots,\tau_q\}.
\end{equation}

In order to prove (\ref{eq26}), let us examine products of $\sigma$s and $\tau$s of length $2$.
Clearly, $\sigma_u^{\LEFT}\sigma_v^{\LEFT}=\sigma_u^{\LEFT}$ for any $u,v\in \{2,\ldots,m\}$.
Next note that $\im(\sigma_u^{\RIGHT})=\{1,\ldots,p+1\}$, and that $\sigma_v^{\RIGHT}$ acts identically on it.
Hence $\sigma_u^{\RIGHT}\sigma_v^{\RIGHT}=\sigma_u^{\RIGHT}$, and we conclude that
\begin{equation}
\label{eq26a}
\sigma_u\sigma_v=\sigma_u.
\end{equation}

Let us now examine the product $\tau_u\tau_v$ ($1\leq u,v\leq q$).
If $u=v$ then of course $\tau_u\tau_v=\tau_u$; so let us suppose $u\neq v$.
Since $\im(\tau_v^{\LEFT})=\{2v,2v+1\}$, and both these points are mapped to $2u$ by $\tau_u^{\LEFT}$,
it follows that $\tau_u^{\LEFT}\tau_v^{\LEFT}=\rho_{2u}$, the constant mapping with value $2u$.
A similar argument shows that $\tau_u^{\RIGHT}\tau_v^{\RIGHT}=\rho_1$, the constant $1$.
Hence
\begin{equation}
\label{eq26b}
\tau_u\tau_v=\left\{ \begin{array}{ll}
       \tau_u & \mbox{if } u=v\\
       \rho_{2u,1} & \mbox{if } u\neq v.
      \end{array}\right.
\end{equation}

A very similar argument shows that
\begin{equation}
\label{eq26c}
\sigma_u\tau_v=\rho_{u,1}.
\end{equation}
Finally, let us consider the product $\tau_u\sigma_v$.
On the left we have
$$
\tau_u^{\LEFT}\sigma_v^{\LEFT}=
\left\{ \begin{array}{ll}
             \rho_{2u} & \mbox{if } v\neq 2u+1\\
             \mu       & \mbox{if } v=2u+1,
        \end{array}\right.
$$
where
$$
\mu(x)=\left\{ \begin{array}{ll}
                   2u & \mbox{if } x=1\\
                   2u+1 & \mbox{if } x\neq 1.
               \end{array}\right.
$$
Note that $\mu$ is not an idempotent, since
$$
\mu(\mu(1))=2u+1\neq 2u=\mu(1).
$$
On the right $\tau_u^{\RIGHT}\sigma_v^{\RIGHT}$ is either constant $\rho_1$, or else
it maps $\{ p+2u,p+2u+1\}$ onto some $r\in\{2,\ldots,p+1\}$, and everything else, including $r$, to $1$, in which case it is not an idempotent.
So we conclude that 
\begin{equation}
\label{eq26d}
\tau_u\sigma_v=\rho_{2u,1}\in R_{m,n} \mbox{ or }
\tau_u\sigma_v\not\in E(S).
\end{equation}

From (\ref{eq26a})--(\ref{eq26d}) we see that the only products of length $2$ of generators of $S$ that are not in $R_{m,n}\cup\{\sigma_2,\ldots,\sigma_m,\tau_1,\ldots,\tau_q\}$ are some of $\tau_u\sigma_v$, and that these products are not idempotent. It is now straightforward to see that no product of length greater than 2 is going to produce any further new elements of $S$, and hence that (\ref{eq26}) indeed holds.
This completes the proof of \ref{Rel3}, and we can conclude that
$$
\MaxSbgp(\IG(E(S)),\rho_{11})\cong G,
$$
proving Theorem \ref{thm2}.

Theorem \ref{thm1} is proved in exactly the same way. We just omit the assumption that $G$ be finitely presented, at the expense of allowing the presentation (\ref{eq14}) to become infinite. This in turn makes the index sets $I$ and $J$, the family of idempotents $\sigma_u$, $\tau_v$, and ultimately the semigroup $S$, infinite.
Alternatively, Theorem \ref{thm1} can be deduced as a corollary of Theorem \ref{thm3}, which will be proved in the next section.

\section{Second construction}
\label{sec4}

The semigroup $S$ constructed in Section \ref{sec3} is not regular:
the non-constant products $\tau_u\sigma_v$ are easily seen to be non-regular.
Furthermore, there does not exist a regular semigroup $S^\prime$ with the same set of idempotents as $S$.
Indeed, for $\tau_u$ and $\sigma_v$ as above, the sandwich set $S(\tau_u,\sigma_v)$ is easily seen to be empty, implying that $E(S)$ is not a regular biordered set; 
see \cite{nambooripad79} or \cite[Proposition 2.5.1]{howie95}.
So, in order to prove Theorems \ref{thm3} and \ref{thm4} we introduce a new construction.

This construction has two advantages over that described in the previous section:
it is in a sense more compact, and it always yields a regular semigroup.
Its main disadvantage is that the constructed semigroup is finite if and only if the input group is finite.

Let $G$ be an arbitrary group, let $N$ be its (possibly infinite) order, and let $n=N^2$.
We will work inside the semigroup $B_{3,n}=T_3^{\op}\times T_n$, which has the $3\times n$ rectangular band $R_{3,n}$ as its minimal ideal.
The corresponding generators for the maximal subgroup containing $\rho_{11}$, as introduced in Section \ref{sec2}, are
\begin{equation}
\label{eq27aa}
f_{ij}\ (i=1,2,3; j=1,\ldots,n)
\end{equation}
with
\begin{equation}
\label{eq27ab}
f_{1j}=f_{i1}=1\ (i=1,2,3; j=1,\ldots,n).
\end{equation}

The auxiliary matrix this time is
$$
Y=(y_{ij})_{3\times n}=\left( \begin{array}{ccccc}
                                1&1&1&\ldots&1\\
                                1&y_{22}&y_{23}&\ldots&y_{2n}\\
                                1&y_{32}&y_{33}&\ldots&y_{3n}
                        \end{array}\right).
$$
Its entries are the elements of $G$, arranged arbitrarily subject to the condition that every possible column appears (once and only once):
\begin{equation}
\label{eq27ac}
\{ (1,y_{2j},y_{3j})\::\: j=1,\ldots,n\}=\{(1,g,h)\::\: g,h\in G\}.
\end{equation}
In fact, in what follows we will sometimes identify the index set $J=\{1,\ldots,n\}$
and the set $\{(1,g,h)\::\: g,h\in G\}$ of all columns of $Y$.
When we do want to distinguish between the two sets we will write $Y_j$ for the $j$th column:
$Y_j=(1,y_{2j},y_{3j})$. The index set $I$ is, of course, just $\{1,2,3\}$.

Now we proceed to define our extra idempotents.
There are six of them
\begin{equation}
\label{eq27ae}
\sigma_u=(\sigma_u^{\LEFT},\sigma_u^{\RIGHT})\ (u=1,\ldots,6),
\end{equation}
and are given by
\begin{alignat*}{5}
&\sigma_1^{\LEFT}&&=\left({\setlength{\arraycolsep}{1mm}\begin{array}{ccc}1&2&3\\ 1&2&2\end{array}}\right) &&\quad\quad
&&\sigma_1^{\RIGHT}&&:(1,g,h)\mapsto (1,g,g)
\\
&\sigma_2^{\LEFT}&&=\left({\setlength{\arraycolsep}{1mm}\begin{array}{ccc}1&2&3\\ 1&2&1\end{array}}\right) &&\quad 
&&\sigma_2^{\RIGHT}&&:(1,g,h)\mapsto (1,g,1)
\\
&\sigma_3^{\LEFT}&&=\left({\setlength{\arraycolsep}{1mm}\begin{array}{ccc}1&2&3\\ 1&1&3\end{array}}\right) && \quad 
&&\sigma_3^{\RIGHT}&&:(1,g,h)\mapsto (1,1,h)
\\
&\sigma_4^{\LEFT}&&=\left({\setlength{\arraycolsep}{1mm}\begin{array}{ccc}1&2&3\\ 1&3&3\end{array}}\right) && \quad 
&&\sigma_4^{\RIGHT}&&:(1,g,h)\mapsto (1,h,h)
\\
&\sigma_5^{\LEFT}&&=\left({\setlength{\arraycolsep}{1mm}\begin{array}{ccc}1&2&3\\ 2&2&3\end{array}}\right) && \quad 
&&\sigma_5^{\RIGHT}&&:(1,g,h)\mapsto (1,1,hg^{-1})
\\
&\sigma_6^{\LEFT}&&=\left({\setlength{\arraycolsep}{1mm}\begin{array}{ccc}1&2&3\\ 3&2&3\end{array}}\right) && \quad 
&&\sigma_6^{\RIGHT}&&:(1,g,h)\mapsto (1,gh^{-1},1).
\end{alignat*}
Note that only $\sigma_5^{\RIGHT}$ and $\sigma_6^{\RIGHT}$ actually depend on the group $G$.

Let us illustrate this
with an example.
We take 
$$G=K_4=\langle a,b\: |\: a^2=b^2=1,\ ab=ba\rangle=\{1,a,b,c\},$$
the Klein four-group.
The index sets are $I=\{1,2,3\}$ and $J=\{1,\ldots,16\}$, and the auxiliary matrix is
$$
Y=\left(\begin{array}{cccccccccccccccc}
1&1&1&1&1&1&1&1&1&1&1&1&1&1&1&1\\
1&1&1&1&a&a&a&a&b&b&b&b&c&c&c&c\\
1&a&b&c&1&a&b&c&1&a&b&c&1&a&b&c
\end{array}\right).
$$
If we want to compute, say, $(12)\sigma_6^{\RIGHT}$ we proceed as follows.
Take the $12$th column of $Y$ -- that is $(1,b,c)$.
Transforming it via $(1,g,h)\mapsto (1,gh^{-1},1)$ yields $(1,a,1)$, which is column $5$.
Thus $(12)\sigma_6^{\RIGHT}=5$.

The full list of all the mappings is:

\begin{alignat*}{5}
&\sigma_1^{\LEFT}&=\left({\setlength{\arraycolsep}{1mm}\begin{array}{ccc}1&2&3\\ 1&2&2\end{array}}\right)&\quad&
\sigma_1^{\RIGHT}&=\left({\setlength{\arraycolsep}{1mm}\begin{array}{cccccccccccccccc} 
                    1&2&3&4&5&6&7&8&9&10&11&12&13&14&15&16\\ 
                    1&1&1&1&6&6&6&6&11&11&11&11&16&16&16&16
               \end{array}}\right)
\\
&\sigma_2^{\LEFT}&=\left({\setlength{\arraycolsep}{1mm}\begin{array}{ccc}1&2&3\\ 1&2&1\end{array}}\right)&&
\sigma_2^{\RIGHT}&=\left({\setlength{\arraycolsep}{1mm}\begin{array}{cccccccccccccccc}
                    1&2&3&4&5&6&7&8&9&10&11&12&13&14&15&16\\ 
                    1&1&1&1&5&5&5&5&9&9&9&9&13&13&13&13
               \end{array}}\right)
\\
&\sigma_3^{\LEFT}&=\left({\setlength{\arraycolsep}{1mm}\begin{array}{ccc}1&2&3\\ 1&1&3\end{array}}\right)&&
\sigma_3^{\RIGHT}&=\left({\setlength{\arraycolsep}{1mm}\begin{array}{cccccccccccccccc}
                    1&2&3&4&5&6&7&8&9&10&11&12&13&14&15&16\\ 
                    1&2&3&4&1&2&3&4&1&2&3&4&1&2&3&4
               \end{array}}\right)
\\
&\sigma_4^{\LEFT}&=\left({\setlength{\arraycolsep}{1mm}\begin{array}{ccc}1&2&3\\ 1&3&3\end{array}}\right)&&
\sigma_4^{\RIGHT}&=\left({\setlength{\arraycolsep}{1mm}\begin{array}{cccccccccccccccc}
                    1&2&3&4&5&6&7&8&9&10&11&12&13&14&15&16\\ 
                    1&6&11&16&1&6&11&16&1&6&11&16&1&6&11&16
               \end{array}}\right)
\\
&\sigma_5^{\LEFT}&=\left({\setlength{\arraycolsep}{1mm}\begin{array}{ccc}1&2&3\\ 2&2&3\end{array}}\right)&&
\sigma_5^{\RIGHT}&=\left({\setlength{\arraycolsep}{1mm}\begin{array}{cccccccccccccccc}
                    1&2&3&4&5&6&7&8&9&10&11&12&13&14&15&16\\ 
                    1&2&3&4&2&1&4&3&3&4&1&2&4&3&2&1
               \end{array}}\right)
\\
&\sigma_6^{\LEFT}&=\left({\setlength{\arraycolsep}{1mm}\begin{array}{ccc}1&2&3\\ 3&2&3\end{array}}\right)&&
\sigma_6^{\RIGHT}&=\left({\setlength{\arraycolsep}{1mm}\begin{array}{cccccccccccccccc}
                    1&2&3&4&5&6&7&8&9&10&11&12&13&14&15&16\\ 
                    1&5&9&13&5&1&13&9&9&13&1&5&13&9&5&1
               \end{array}}\right).
\end{alignat*}

If we were to change the group $G$, say to $G=C_4=\langle a\:|\: a^4=1\rangle=\{1,a,a^2,a^3\}=\{1,a,b,c\}$, the cyclic group of order $4$, only the mappings $\sigma_5^{\RIGHT}$ and $\sigma_6^{\RIGHT}$ would change:
\begin{alignat*}{2}
&\sigma_5^{\RIGHT}&&=\left({\setlength{\arraycolsep}{1mm}\begin{array}{cccccccccccccccc}
                    1&2&3&4&5&6&7&8&9&10&11&12&13&14&15&16\\ 
                    1&2&3&4&4&1&2&3&3&4&1&2&2&3&4&1
               \end{array}}\right)
\\
&\sigma_6^{\RIGHT}&&=\left({\setlength{\arraycolsep}{1mm}\begin{array}{cccccccccccccccc}
                    1&2&3&4&5&6&7&8&9&10&11&12&13&14&15&16\\ 
                    1&13&9&5&5&1&13&9&9&5&1&13&13&9&5&1
               \end{array}}\right).
\end{alignat*}

We now return to the main argument. A routine verification shows that the semigroup generated by
$\{\sigma_u\::\: u=1,\ldots,6\}$ has $21$ elements:
$\sigma_u$ ($u=1,\ldots,6$) themselves, three elements $\rho_{11}$, $\rho_{21}$, $\rho_{31}$ belonging to the minimal ideal $R_{3,n}$, and twelve further elements:
\begin{alignat*}{3}
&\sigma_7& &=\sigma_1\sigma_3 &&=( \left({\setlength{\arraycolsep}{1mm}\begin{array}{ccc} 1&2&3\\ 1&1&2\end{array}}\right),
                             (1,g,h)\mapsto  (1,1,g) )
\\
&\sigma_8& &=\sigma_2\sigma_5 & &=( \left({\setlength{\arraycolsep}{1mm}\begin{array}{ccc} 1&2&3\\ 2&2&1\end{array}}\right),
                             (1,g,h)\mapsto  (1,1,g^{-1}) )
\\
&\sigma_9&&=\sigma_3\sigma_6&&=( \left({\setlength{\arraycolsep}{1mm}\begin{array}{ccc} 1&2&3\\ 3&1&3\end{array}}\right),
                             (1,g,h)\mapsto  (1,h^{-1},1) )
\\
&\sigma_{10}&&=\sigma_4\sigma_2&&=( \left({\setlength{\arraycolsep}{1mm}\begin{array}{ccc} 1&2&3\\1&3&1\end{array}}\right),
                             (1,g,h)\mapsto  (1,h,1) )
\\
&\sigma_{11}&&=\sigma_5\sigma_4&&=( \left({\setlength{\arraycolsep}{1mm}\begin{array}{ccc} 1&2&3\\ 2&3&3\end{array}}\right),
                             (1,g,h)\mapsto  (1,hg^{-1},hg^{-1}) )
\\
&\sigma_{12}&&=\sigma_6\sigma_1&&=( \left({\setlength{\arraycolsep}{1mm}\begin{array}{ccc} 1&2&3\\ 3&2&2\end{array}}\right),
                             (1,g,h)\mapsto  (1,gh^{-1},gh^{-1}) )
\\
&\sigma_{13}&&=\sigma_1\sigma_3\sigma_6&&=( \left({\setlength{\arraycolsep}{1mm}\begin{array}{ccc} 1&2&3\\ 2&1&2\end{array}}\right),
                             (1,g,h)\mapsto  (1,g^{-1},1) )
\\
&\sigma_{14}&&=\sigma_2\sigma_5\sigma_4&&=( \left({\setlength{\arraycolsep}{1mm}\begin{array}{ccc} 1&2&3\\ 2&1&1\end{array}}\right),
                             (1,g,h)\mapsto  (1,g^{-1},g^{-1}) )
\\
&\sigma_{15}&&=\sigma_3\sigma_6\sigma_1&&=( \left({\setlength{\arraycolsep}{1mm}\begin{array}{ccc} 1&2&3\\ 3&1&1\end{array}}\right),
                             (1,g,h)\mapsto  (1,h^{-1},h^{-1}) )
\\
&\sigma_{16}&&=\sigma_4\sigma_2\sigma_5&&=( \left({\setlength{\arraycolsep}{1mm}\begin{array}{ccc} 1&2&3\\ 3&3&1\end{array}}\right),
                             (1,g,h)\mapsto  (1,1,h^{-1}) )
\\
&\sigma_{17}&&=\sigma_5\sigma_4\sigma_2&&=( \left({\setlength{\arraycolsep}{1mm}\begin{array}{ccc} 1&2&3\\ 2&3&2\end{array}}\right),
                             (1,g,h)\mapsto  (1,hg^{-1},1) )
\\
&\sigma_{18}&&=\sigma_6\sigma_1\sigma_3&&=( \left({\setlength{\arraycolsep}{1mm}\begin{array}{ccc} 1&2&3\\ 3&3&2\end{array}}\right),
                             (1,g,h)\mapsto  (1,1,gh^{-1}) ).
\end{alignat*}
The elements $\{ \sigma_1,\ldots,\sigma_{18}\}$ form a single $\GD$-class $D$ illustrated in Figure \ref{fig2}.

\begin{figure}
\begin{center}
\begin{picture}(45,45)
\thicklines

\put(0,0){\line(1,0){45}}
\put(0,15){\line(1,0){45}}
\put(0,30){\line(1,0){45}}
\put(0,45){\line(1,0){45}}
\put(0,0){\line(0,1){45}}
\put(15,0){\line(0,1){45}}
\put(30,0){\line(0,1){45}}
\put(45,0){\line(0,1){45}}

\put(3,38){\makebox(3,3){$\sigma_1$}}
\put(4.5,39.5){\circle{6}}
\put(8,32){\makebox(3,3){$\sigma_{14}$}}
\put(18,38){\makebox(3,3){$\sigma_2$}}
\put(19.5,39.5){\circle{6}}
\put(23,32){\makebox(3,3){$\sigma_{13}$}}
\put(33,38){\makebox(3,3){$\sigma_7$}}
\put(38,32){\makebox(3,3){$\sigma_{8}$}}
\put(3,23){\makebox(3,3){$\sigma_4$}}
\put(4.5,24.5){\circle{6}}
\put(8,17){\makebox(3,3){$\sigma_{15}$}}
\put(18,23){\makebox(3,3){$\sigma_9$}}
\put(23,17){\makebox(3,3){$\sigma_{10}$}}
\put(33,23){\makebox(3,3){$\sigma_3$}}
\put(34.5,24.5){\circle{6}}
\put(38,17){\makebox(3,3){$\sigma_{16}$}}
\put(3,8){\makebox(3,3){$\sigma_{11}$}}
\put(8,2){\makebox(3,3){$\sigma_{12}$}}
\put(18,8){\makebox(3,3){$\sigma_6$}}
\put(19.5,9.5){\circle{6}}
\put(23,2){\makebox(3,3){$\sigma_{17}$}}
\put(33,8){\makebox(3,3){$\sigma_5$}}
\put(34.5,9.5){\circle{6}}
\put(38,2){\makebox(3,3){$\sigma_{18}$}}
\end{picture}
\caption{The $\GD$-class $D$, with the idempotents circled.}
\label{fig2}
\end{center}
\end{figure}

It follows that the semigroup
$$
S=\langle R_{3,n}\cup \{\sigma_1,\ldots,\sigma_6\}\rangle=R_{3,n}\cup D\leq B_{3,n}
$$
has the following properties:
\begin{itemize}
\item
$S$ is regular;
\item
$S$ has two $\GD$-classes: $R_{3,n}$ and $D$;
\item
$S$ has precisely six idempotents ($\sigma_1,\ldots,\sigma_6$) outside $R_{3,n}$;
\item\
$S$ is finite if and only if $R_{3,n}$ is finite, which is the case if and only if $G$ is finite.
\end{itemize}

We now follow the process outlined in Section \ref{sec2} to prove that the maximal subgroup of $\IG(E(S))$ containing $\rho_{11}\in R_{3,n}$ is isomorphic to $G$.
So far we have got the generators (\ref{eq27aa}) and relations (\ref{eq27ab}).
Any further relations arise from the squares singularised by $\sigma_u$, $u=1,\ldots,6$.

Since $\im(\sigma_1^{\LEFT})=\{1,2\}$ it follows that the left-right squares singularised by $\sigma_1$ are of the form $(1,2;j,l)$, where $Y_j=(1,g,g)$, $Y_l=(1,g,h)$ for some $g,h\in G$.
These are flush top squares yielding the relations
\begin{equation}
\label{eq27ba}
f_{2j}=f_{2l}\ (\text{whenever } y_{2j}=y_{2l}).
\end{equation}
The left-right squares singularised by $\sigma_2$ yield exactly the same relations, while
those singularised by $\sigma_3$ and $\sigma_4$ yield 
\begin{equation}
\label{eq27bb}
f_{3j}=f_{3l}\ (\text{whenever } y_{3j}=y_{3l}).
\end{equation}
Next note that $\im(\sigma_1^{\RIGHT})=\{(1,g,g)\::\: g\in G\}$; so the up-down squares singularised by $\sigma_1$ have the form $(2,3; j,l)$, where $Y_j=(1,g,g)$, $Y_l=(1,h,h)$.
For $j=1$ we obtain the flush left square $(2,3;1,l)$, yielding the relation
\begin{equation}
\label{eq27bc}
f_{2l}=f_{3l}\ (\text{whenever } y_{2l}=y_{3l}).
\end{equation}
The relation $f_{2j}^{-1}f_{2l}=f_{3j}^{-1}f_{3l}$ produced from a general square $(2,3;j,l)$ is an easy consequence of (\ref{eq27bc}).

Combining (\ref{eq27ba}), (\ref{eq27bb}), (\ref{eq27bc}) together gives us precisely the relations
\begin{equation}
\label{eq27d}
f_{ij}=f_{kl}\ (\text{whenever } y_{ij}=y_{kl}),
\end{equation}
i.e. we have completed step \ref{Rel1}, and can identify each generator $f_{ij}$ with the entry $y_{ij}$ of the auxiliary matrix (considered as a formal symbol).

The up-down squares singularised by $\sigma_2$, $\sigma_3$, $\sigma_4$ do not yield any relations over and above (\ref{eq27d}). So there remains to analyse the relations yielded by $\sigma_5$ and $\sigma_6$.
From $\im(\sigma_5^{\LEFT})=\{2,3\}$ it follows that the left-right squares singularised by
$\sigma_5$ are of the form $(2,3; j,l)$, where $Y_j=(1,1,hg^{-1})$ and $Y_l=(1,g,h)$.
This is a $3/4$ square yielding the relation $(hg^{-1})\cdot g=h$.
Clearly, as $g$ and $h$ range through $G$, we obtain the Cayley table of $G$ (considered as a presentation), accomplishing \ref{Rel2}.
The up-down squares for $\sigma_5$ have the form $(2,1; j,l)$, where $Z_j=(1,1,g)$ and $Z_l=(1,1,h)$,
and yield the trivial relation $1=1$. Likewise, $\sigma_6$ yields no further relation.
Since $S$ has no further idempotents, we have accomplished \ref{Rel3}. This completes the proof
of Theorems \ref{thm3} and \ref{thm4}.

\section{Concluding remarks}
\label{sec5}

During the work on this project we have implemented in \textsf{GAP} \cite{GAP4} the 
Rei\-de\-meis\-ter--Schreier type rewriting process described in Section \ref{sec_prelims} and used in Section \ref{sec_IGpres}. 
This has enabled us to gather considerable `experimental data' and test several early conjectures. 
The output from any Reidemeister--Schreier type rewriting has a large number of generators and defining relations.
Thus, the Tietze Transformations programme, which is a part of the standard \textsf{GAP} distribution,
and which is in the \textsf{GAP} manual credited back to the work of 
Havas, Robertson et al. \cite{havas69,havas84,robertson88},
proved an invaluable tool.

It is well known that there exists a finitely presented group 
with an unsolvable word problem.
Combining such a group with Theorem \ref{thm2} yields:

\begin{corollary}
There exists a free idempotent generated semigroup $F$ arising from a finite semigroup
such that the word problem for $F$ is unsolvable.
\end{corollary}

Such a free idempotent generated semigroup $F$ would be non residually finite and non-automatic as well.

The main open question remaining, as mentioned in the Introduction, is:

\begin{problem}
Is it true that every finitely presented group is a maximal subgroup of some free regular idempotent generated semigroup arising from a finite semigroup?
\end{problem}

All our examples were obtained by first fixing a rectangular band $R$, and then constructing our semigroup $S$ as an ideal extension of $R$.
We remark that for any fixed finite $R$ we can obtain only finitely many maximal subgroups in free idempotent semigroups arising from ideal extensions from $R$. Indeed, $R$ determines the generators of the maximal subgroups, and the relators arise from singular squares, and there are only finitely many of them. Thus we are lead to ask:

\begin{problem}
Given a finite $0$-simple semigroup $R$ and and idempotent $e\in R$, describe all the groups which arise as maximal subgroups containing $e$ of the free idempotent semigroups $\IG(S)$ where $S$ is a finite ideal extension of $R$.
\end{problem}

\begin{acknowledgement}
The authors would like to thank an anonymous referee,
whose comments led them to state and prove Theorem \ref{lismaintheorem} in its present form.
We would also like to thank Stuart W. Margolis for insightful correspondence and helpful suggestions.
\end{acknowledgement}

\begin{flushleft}
Centro de \'{A}lgebra da Universidade de Lisboa \\ 
Av. Prof. Gama Pinto 2  \\
1649-003 Lisboa,  Portugal.\\
\smallskip
\texttt{rdgray@fc.ul.pt}
\end{flushleft}

\begin{flushleft}
School of Mathematics and Statistics\\
University of St Andrews\\
St Andrews KY16 9SS\\
Scotland, U.K.\\
\smallskip
\texttt{nik@mcs.st-and.ac.uk}
\end{flushleft}

\end{document}